\documentclass[a4paper,12pt]{article}
\usepackage{amsmath}
\usepackage{amssymb}
\usepackage{tabularx}
\usepackage{enumerate}

\newtheorem{thm}{Theorem}[section]
\newtheorem{lem}[thm]{Lemma}
\newtheorem{cor}[thm]{Corollary}
\newtheorem{prop}[thm]{Proposition}
\newtheorem{rmk}{Remark}[section]

\newtheorem{pppp}{Proof}

\newcommand{\qed}{\hspace{1em}\mbox{\raisebox{0.65ex}{\fbox{}}}}

\numberwithin{equation}{section}

\newcommand{\be}{\begin{equation}}
\newcommand{\ee}{\end{equation}}
\newcommand\bes{\begin{eqnarray}} \newcommand\ees{\end{eqnarray}}
\newcommand{\bess}{\begin{eqnarray*}}
\newcommand{\eess}{\end{eqnarray*}}

\newcommand{\bpf}{{\bf Proof:\ \ }}
\newcommand{\epf}{\mbox{}\hfill $\Box$}

\begin{document}

\thispagestyle{empty}

\title{The spreading front of invasive species in favorable habitat or unfavorable
habitat\thanks{The work is supported by the PRC grant NSFC 11271197, NSFC 11371311 and NSF of the Higher Education Institutions of Jiangsu Province
(12KJD110008).}}
\date{\empty}

\author{Chengxia Lei$^{a}$, Zhigui Lin$^b$ \thanks{Corresponding author. Email: zglin68@hotmail.com} and Qunying Zhang$^{b}$\\
{\small $^a$Department of Mathematics, University of Science and Technology of China,}\\
{\small Hefei 230026, China}\\
{\small $^b$School of Mathematical Science, Yangzhou University, Yangzhou 225002, China}
}
 \maketitle

\begin{quote}
\noindent
{\bf Abstract.} { 
\small Spatial heterogeneity and habitat characteristic are shown to determine the asymptotic
profile of the solution to a reaction-diffusion model with free boundary,
which describes the moving front of the invasive species.
A threshold value $R_0^{Fr}(D,t)$ is introduced to determine the spreading and vanishing of the invasive species.
We prove that if $R_0^{Fr}(D,t_0)\geq 1$ for some $t_0\geq 0$, the spreading must happen; while if $R_0^{Fr}(D,0)<1$,
the spreading is also possible. Our results show that the species in the favorable habitat can establish itself
if the diffusion is slow or the occupying habitat is large.
In an unfavorable habitat, the species dies out if the initial value of the species is small.
However, big initial number of the species is benefit for the species to survive.
When the species spreads in the whole habitat, the asymptotic spreading speed is given.
Some implications of these theoretical results are also discussed. }



\end{quote}

\section{Introduction}

There have been some recent studies on reaction-diffusion models to understand
the nature of spreading of the invasive species. The spreading of species
from their native habitats to alien environments is a serious
threat to biological diversity \cite{MCH}. Subsequently,
many mathematical models have been constructed to investigate how
raw species survive in the habitat \cite{SK}. Among those models, there was a
well-known model, which is described by the diffusive logistic equation over the
entire space $\mathbb{R}^n$:
\begin{equation}
u_{t}-d \Delta u=u(a-b u),\; x\in {\mathbb R}^n,\; t>0, \label{fish}
\end{equation}
where $u=u(x,t)$ is the population density of an invasive species
 with diffusion rate $d$, intrinsic growth rate $a$ and habitat carrying capacity $a/b$.
For space dimension $n = 1$,
traveling wave solutions have been found by Fisher \cite{F} and Kolmogorov et al \cite{KPP}. That is, for any $c\geq c^*:= 2\sqrt{ad}$, there exists
a solution $u(x,t) := W(x-ct)$ with the property that
$$W'(y)<0,\ y\in {\mathbb R}^1,\ W(-\infty)=a/b,\ W(\infty)=0;$$
no such solution exists if $c < c^*$.
The number $c^*$ is then called the minimal speed of the
traveling waves. The related research and recent developments can be found, for example, in \cite{DL2, HH} and the
references therein.

However, the solution to problem (\ref{fish}) with any nontrivial initial population $u(x,0)$ is always positive everywhere,
which means that any invasive species can establish itself in any new environment. This contrasts sharply with
numerous empirical evidences, for example, let us see a biological control programme for broom began in New Zealand in 1981.
A field experiment was used to manipulate the critical first stages of the invasion
of the psyllid, Arytainilla spartiophila. Fifty-five releases were made along a linear transect 135 km long, six years after their original
release, psyllids were present in 22 of the 55 release sites \cite{CC}.
The field experiment showed that not all releases of psyllids survived and successful establishment is a complex process.

To describe precisely the spreading front of invasive species,
 Du and Lin \cite{DL} studied the following free boundary problem,
\begin{eqnarray}
\left\{
\begin{array}{lll}
u_{t}-d u_{xx}=u(a-b u),\; &0<x<h(t),\; t>0,   \\
u_x(0,t)=u(h(t),t)=0,& t>0,\\
h'(t)=-\mu u_{x}(h(t),t),\;  & t>0, \\
u(x,0)=u_{0}(x), & 0\leq x\leq h_0,
\end{array} \right.
\label{a1}
\end{eqnarray}
where $x = h(t)$ is the free boundary to be determined, $d, a, h_0, \mu$ and $b$ are given
positive constants, the unknown $u(x,t)$ stands for the population density of an invasive species over a one-dimensional
habitat, and the initial function $u_0(x)$ stands for the population of the
species in the early stage of its introduction. It is assumed that the spreading front expands at a speed
that is proportional to the population gradient at the front, which gives rise to the classical
Stefan condition $h'(t) = -\mu u_x(h(t),t)$, the positive constant $\mu$ measures the ability of the invasive species
to transmit and diffuse in the new habitat, see \cite{LIN} in details.

A spreading-vanishing dichotomy was first presented in \cite{DL} for problem (\ref{a1}), namely, as time approaches to infinity, the population
$u(x,t)$ either successfully establishes itself in the new environment (called spreading),
in the sense that $h(t)\to \infty$ and $u(x,t)\to a/b$, or the population fails to establish and
vanishes eventually (called vanishing), namely $h(t)\to h_\infty\leq \frac \pi 2\sqrt{\frac d a}$ and $u(x,t)\to 0$.

It was also shown that if spreading occurs, for large time, the spreading speed approaches
a positive constant $k_0$, i.e., $h(t) = [k_0 + o(1)]t$ as $t\to \infty$. $k_0$ is then
called the asymptotic spreading speed, which is uniquely determined by an auxiliary
elliptic problem induced from (\ref{a1}). Furthermore, they found that $k_0<c^*$,
where $c^*(:=2\sqrt{ad})$ is the minimal speed of the traveling waves.
Hereafter, Du and Guo \cite{DG} extended the free boundary problem (1.2) to a higher dimension domain.

The classical example of biological invasion is the spreading of the muskrat in Europe \cite{LHM}.
As we know that the original habitat of the muskrat was North America before 1905,
and there was no one in Europe.  A few muskrats have been brought to Prague until 1905
and five of them were escaped from a farm.
With the reproduction of the muskrats, today Europe is believed to contain more than 100,000,000 muskrats.
Skellam calculated the area of the muskrat range according to a map obtained from field data
and found that the spreading radius eventually exhibits a linear growth curve against time \cite{Ske}.

Another successful invasion is the spreading of the cane toad in Australia.
To help cane-growers in Queensland control beetles and increase yields,
101 young Hawaiian toads were brought to Australia in August 1935.
After that they have rapidly multiplied in population and have steadily expanded their range.
It is estimated that the number is  over 200 million and toad migrates at an average
of 40 kilometers per year \cite{Ty}.

Since the work of Du and Lin \cite{DL}, there have been many
theoretical developments on the free boundary problem in homogeneous
environment. For example,  Kaneko and Yamada \cite{KY} considered a free boundary problem for a general reaction-diffusion
equation with Dirichlet conditions on both fixed and free boundaries. Du and Lou \cite{DL3} discussed a two free boundaries problem with
 a general nonlinear term. In \cite{GLL}, Gu, Lin and Lou studied how advection
term $(\beta u_x)$ affects the asymptotic spreading speeds when
spreading occurs $(0<\beta<2\sqrt{d})$. See also \cite{DGP} for diffusive logistic model in time-periodic
environment, \cite{PZ} for diffusive logistic model with seasonal succession,
\cite{LLW} for information diffusion in online social networks,
\cite{Wa, WZ} for Lotka-Volterra type prey-predator model and
\cite{DL1, GW} for Lotka-Volterra type competition model.

But these models are not exactly describe the survival of species in real environment, for example,
 some part of the habitat has been polluted or destroyed.
To illustrate how the spatial heterogeneity affect persistence or extinction of a species,
Cantrell and Cosner \cite{CC1} proposed the following diffusive logistic equation in a fixed domain $\Omega$:
 \begin{eqnarray}
\left\{
\begin{array}{lll}
u_{t}-d \Delta u=u(m(x)-cu),\; &(x,t)\in\Omega\times(0,+\infty),   \\
u(x,t)=0,&(x,t)\in\partial\Omega\times(0,+\infty),\\
u(x,0)=u_{0}(x), & x\in \Omega,
\end{array} \right.
\label{a2}
\end{eqnarray}
where $u(x,t)$ represents the density of the species, $c,d$ are positive constants and $m(x)\in L^\infty(\Omega)$,
and
$$\Omega_m:=\left\{x\in\Omega: m(x)\geq0\right\}$$
is nonempty and not equal to the whole domain $\Omega$.
For a fixed diffusion rate $d$, they concluded that there has a ``ancellation"
effect if the favourable and unfavourable habitats were closely intermingled,
and a small number of large favourable habitats is better for the species to survive than many small ones.

Inspired by the former work, we will focus on the impact of spatial feature of environment on the spreading and vanishing
of an invasive species with a free boundary describing the moving front.
 For simplicity, we assume the environment is radially symmetric and investigate the behavior of the positive solution
$(u(r, t); h(t))$ with $r(:=|x|, x\in\mathbb{R}^n)$ to the following problem
\begin{eqnarray}
\left\{
\begin{array}{lll}
u_{t}-D \Delta u=u(b(r) -d(r)-\beta (r)u),\; &0<r<h(t),\; t>0,   \\
u_r(0,t)=u(r,t)=0,&r=h(t),\, t>0,\\
h'(t)=-\mu u_{r}(h(t),t),\;   & t>0, \\
u(r,0)=u_{0}(r), & 0\leq r\leq h_0,
\end{array} \right.
\label{a3}
\end{eqnarray}
where $\triangle u=u_{rr}+\frac{n-1}{r}u_r$, $r=h(t)$ is the moving
boundary to be defined,  $h_0,\ D$ and $\mu $ are positive constants as above, $b(r)$, $d(r)$
and $\beta(r)$ are positive H$\ddot{o}$lder continuous functions which account for the birth rate, death rate and crowding strength
of the species at $r$, respectively.
In the paper, we assume that there exist positive constants $b_1$ and $b_2$ such that $b_1\leq b(r),d(r),\beta(r)\leq b_2$
for $r\in[0,+\infty)$. The initial function $u_0(r)$ is nonnegative and satisfies
\begin{equation}
u_0\in C^2([0, h_0]),\, u'_0(0)=u_0(h_0)=0\,  \textrm{and} \ u_0(r)>0,\, r\in [0, h_0),
\label{Ae2}
\end{equation}
where the condition (\ref{Ae2}) indicates that in the early stage of its introduction,
the species exists in the area with $r\in [0, h_0)$,
beyond the free boundary $r=h(t)$, there is no invasive species.

To describe the feature of environment, as in \cite{AL, L1},
we say that $r$ is a {\bf favorable site} if the local birth rate $b(r)$ is greater
than the local death rate $d(r)$. An {\bf unfavorable site} is defined in a similar manner.
Denote the favorable set and the unfavorable set as following:
$$F^+=\{r\in(0,+\infty):b(r)>d(r)\} \,  \textrm{and}  \  F^-=\{r\in(0,\infty):b(r)<d(r)\}.$$
The habitat $B_R$ (a ball with radius $R$) is characterized as {\bf favorable}
( or {\bf unfavorable} ) if the spatial average $(\frac{1}{|B_R|}\int_{B_R}b(r)dr)$
of the birth rate is greater than ( or less than ) the spatial
average $(\frac{1}{|B_R|}\int_{B_R}d(r)dr)$ of the death rate, respectively.

To establish the theoretical conclusion,
we introduce the threshold value
$$R_0^{Di}(D, B_R):=\ \sup_{\phi\in H^1_0(B_R),\phi\neq 0}\left \{\frac{\int_{B_R} b(r)
\phi^2dr}{\int_{B_R} (D|\triangledown \phi|^2+d(r)\phi^2)d r}\right\}$$
for the corresponding problem (\ref{B11}) in $B_R$ with null Dirichlet boundary condition, and then
the threshold value
$$R_0^{Fr}(D,t):=R_0^{Di}(D, B_{h(t)})=\ \sup_{\phi\in H^1_0(B_{h(t)}),\phi\neq 0} \left \{\frac{\int_{B_{h(t)}} b(r)
\phi^2dr}{\int_{B_{h(t)}} (D|\triangledown \phi|^2+d(r) \phi^2)dr}\right\}$$
is defined for the free boundary problem (\ref{a3}).

In this paper, we mainly study the asymptotic behavior of problem (\ref{a3}).
When spreading occurs, we consider an auxiliary elliptic problem
\begin{eqnarray}
-D\Delta u=u[b(r)-d(r)-\beta(r)u],\ r\in  [0,+\infty).
\label{b121212}
\end{eqnarray}
If we assume that
$$(H)\;\;\;\;\lim_{r\rightarrow\infty}(b(r)-d(r))=\alpha>0,\;\;\;\;\;\;
\lim_{r\rightarrow\infty}\beta(r)=\beta>0,$$
then problem (\ref{b121212}) admits a unique positive solution $\widetilde{u}(r)$ (\cite{DM}).
Let $(u(r,t),h(t))$ be the solution of problem (\ref{a3}), we say

(D1) {\it Spreading of the invasive species} in the heterogeneous environment if
$$h_\infty=\infty \;\;\mbox{and}\;\; \lim_{t\rightarrow+\infty}u(r,t)=\widetilde{u}(r)$$
locally uniformly in $[0,+\infty)$, where $\widetilde{u}(r)$ is the unique positive solution of problem $(\ref{b121212})$.

(D2) {\it Vanishing of the invasive species} in the heterogeneous environment if
$$h_\infty<\infty\;\;\mbox{and}\;\; \lim_{t\rightarrow+\infty}\|u(\cdot,t)\|_{C([0,h(t)])}=0.$$

The remainder of this paper is organized as follows.
In the next section, the global existence and uniqueness of the solution to problem (\ref{a3})
are proved by contraction mapping theorem, comparison principle is also employed.
Section 3 is devoted to the threshold value $R_0^{Di}(D,B_{R})$ for the null Dirichlet boundary and
$R_0^{Fr}(D,t)$ for the free boundary. Their properties are discussed in Theorems 3.1 and 3.3.
Sufficient conditions for the invasive species to spread are given in Section 4. We demonstrate that the species will spread in a favourable
habitat if the diffusion is slow or the habitat occupancy is large. Section 5 deals with the case $R_0^{Fr}(D,0)<1$
and the spreading-vanishing dichotomy is given. In the case of $R_0^{Fr}(D,0)<1$,
the vanishing happens provided the initial value $u_0(r)$ is sufficiently small;
the spreading happens in this circumstance if the favorable set $F^+$ is nonempty and the initial value is large enough.
The diffusion rate, the initial value and the original habitat play a significant role in determine the spreading-vanishing dichotomy.
In Section 6, we given the estimate of the spreading speed when spreading happens. Finally, a brief discussion is given in Section 7.

\section{Existence and uniqueness}

In this section, we first present the local existence and
uniqueness result by contraction mapping theorem, and then use
suitable estimates to show that the solution is well defined for all $t>0$.

\begin{thm} For any given $u_0$ satisfying $(\ref{Ae2})$,
and any $\alpha \in (0, 1)$, there is a $T_0>0$ such that problem \eqref{a3}
admits a unique solution
$$(u, h)\in C^{ 1+\alpha,(1+\alpha)/2}([0,h(t)]\times [0, T_0])\times C^{1+\alpha/2}([0,T_0]);$$
moreover,
\begin{eqnarray} \|u\|_{C^{1+\alpha,
(1+\alpha)/2
}([0,h(t)]\times [0, T_0])}+||h\|_{C^{1+\alpha/2}([0,T_0])}\leq
C,\label{b12}
\end{eqnarray}
where $C$ and $T_0$ only depend on $h_0, \alpha$ and $\|u_0\|_{C^{2}([0, h_0])}$.
\end{thm}
\bpf We can obtain this result by the same arguments as in \cite{CF, DL} with obvious modification. For briefly,
we omit it here.
\epf

By interior Schauder parabolic estimates in \cite{LSU, L}, we have additional regularity for $(u,h)$,
namely,  $u\in C^{2+\alpha,1+\alpha/2}([0, h(t)]\times(0,T_0])$ and $h\in C^{1+\alpha/2}((0,T_0])$.

Next, we collect some basic facts which will be used to show that the local
solution obtained in Theorem 2.1 can be extended to all $t>0$.

\begin{lem} Let $(u, h)$ be a solution to problem \eqref{a3} defined for $t\in (0,T_0]$
for some $T_0>0$.
Then there exists a constant $C_1$ independent of $T_0$ such
that
\[
0<u(r, t)\leq C_1\; \mbox{ for } 0\leq r<
h(t),\; t\in (0,T_0]. \]
Moreover, there exists $C_2$ dependent of $C_1$ but independent of $T_0$ such that
\[
 0<h'(t)\leq C_2 \; \mbox{ for } \; t\in (0,T_0]. \]
\end{lem}
\bpf
Noting that $u\geq 0$ in $[0, h(t))\times[0, T_0]$ as long as the solution exists.
Using the strong maximum principle (see \cite{PW}) to the equations of $u$ in $(0, h(t))\times
[0, T_0]$, we immediately obtain
 $$u(r, t)>0\;\; \textrm{for} \ 0\leq r<h(t),\, 0< t\leq T_0.$$
Since $u(r,t)$ satisfies
\begin{eqnarray*}
\left\{
\begin{array}{lll}
u_{t}-D\Delta u=( b(r) -d(r))u-\beta(r) u^2,\; &
0<r<h(t),\, 0<t\leq T_0,  \\
u'(0, t)=u(r,t)=0,\quad & r=h(t),\, 0<t\leq T_0,\\
u(r,0)=u_0(r),\; &0\leq r\leq h_0,
\end{array} \right.
\end{eqnarray*}
the maximum principle implies that
$$u\leq C_1:=\max\left\{
\frac{\max_{r\in[0,h(t)]}( b(r) -d(r))}{\min_{r\in[0,h(t)]}\beta(r)},\
||u_0(r)||_{L^\infty ([0, h_0])}\right\}$$ in $[0, h(t)]\times [0, T_0]$.

It follows from the proof of Lemma 2.2 in \cite{DL} that we can define an auxiliary function
$$w(r,t)=C_1[2M(h(t)-r)-M^2(h(t)-r)^2]$$
for some suitable $M:=\max\left\{\frac{1}{h_0},\sqrt{\frac{b_2}{2D}},\frac{4\|u_0\|_{C^1([0,h_0])}}{3C_1}\right\}$
such that $w(r,t)\geq u(r,t)$ holds over the domain
$$\{(r,t):\ h(t)-M^{-1}<r<h(t),\ 0<t\leq T_0\}.$$
We then have
$$h'(t)=-\mu u_x(h(t),t)\leq C_2:=2\mu MC_1.$$
By the free boundary condition and using the Hopf Lemma, we have $h'(t)>0$. The proof is complete.
\epf

\begin{rmk}
Noting that $r=h(t)$ is strictly monotonic increasing, and then
there exists  $h_\infty\in (0, +\infty]$ such that $\lim_{t\to +\infty} \ h(t)=h_\infty$.
\end{rmk}

\begin{thm} The solution of problem \eqref{a3} exists and is
unique for all $t\in (0,\infty)$.
\end{thm}
\bpf
By Theorem 2.1, we know that problem \eqref{a3} admits a unique solution in $[0, T_{\max})$ with the maximal time $T_{\max}\leq +\infty$.
If $T_{\max}<+\infty$, we can use the same argument in the proof of Theorem 2.3 in \cite{DL} to derive a contradiction.
Thus the solution of problem \eqref{a3} exists and is unique for all $t\in (0,\infty)$.
\epf

\begin{rmk}
It follows from the uniqueness of the solution to \eqref{a3} and
some standard compactness arguments that the unique solution $(u, h)$
depends continuously on the parameters appearing in \eqref{a3}. This
fact will be used in the sections below.
\end{rmk}

Next, we exhibit the comparison principle.
\begin{lem} (The Comparison Principle \cite{DL})
  Assume that $\underline h, \overline
h\in C^1([0,\infty))$, $\underline u(r,t)\in C([0, \underline h(t)]\times [0, \infty))\cap
C^{2,1}((0, \underline h(t))\times (0, \infty))$,
$\overline u(r,t)\in C([0, \overline h(t)]\times [0, \infty))\cap
C^{2,1}((0, \overline h(t))\times (0, \infty))$ and
\begin{eqnarray}
\left\{
\begin{array}{lll}
\overline u_{t}-D \Delta \overline u\geq (b(r) -d(r))\overline u-\beta (r) \overline u^2, &0<r<\overline h(t), \ t>0,\\
\underline u_{t}-D\Delta \underline u\leq (b(r) -d(r))\underline u-\beta (r)\underline u^2,\; &0<r<\underline h(t),\ t>0, \\
\overline  u(r, t)=0,\;\overline h'(t)\geq -\mu \overline u_r,\quad & r=\overline h(t),\ t>0,\\
\underline  u(r, t)=0,\;\underline h'(t)\leq -\mu \underline u_r,\quad & r=\underline h(t),\ t>0,\\
\underline u(r, 0)\leq u_{0}(r)\leq \overline u(r, 0),\;  &0\leq r\leq h_0,\\
\underline h(0)\leq h_0\leq \overline h(0).   &
\end{array} \right.\label{upp}
\end{eqnarray}
Then the unique solution $(u, h)$ to the free boundary problem $(\ref{a3})$ satisfies
$$h(t)\leq\overline h(t)\, \textrm{for}\, t\in (0, \infty),\ u(r, t)\leq \overline u(r, t)\, \textrm{for}\, (r, t)\in (0,h(t))\times (0, \infty),$$
$$h(t)\geq\underline h(t)\, \textrm{for}\, t\in (0, \infty),\ u(r, t)\geq \underline u(r, t)\, \textrm{for}\, (r, t)\in (0, \underline h(t))\times (0, \infty).$$
\end{lem}

Usually, $(\overline u,\overline h)$ is called a upper solution to problem (\ref{a3})
and $(\underline u,\underline h)$ is a lower solution.

Noting that $\beta (r) \overline u^2\geq 0$, we then have the following Corollary.
\begin{cor}
Suppose that $\overline h\in C^1([0, \infty))$, $\overline u(r,t)\in C([0, \overline h(t)]\times [0, \infty))\cap
C^{2,1}((0, \overline h(t))\times (0, \infty))$, and
\begin{eqnarray*}
\left\{
\begin{array}{lll}
\overline u_{t}-D \Delta \overline u\geq (b (r) -d(r))\overline u, &0<r<\overline h(t), \ t>0,\\
\overline  u=0,\;\overline h'(t)\geq -\mu \overline u_r,\quad & r=\overline h(t),\ t>0,\\
u_{0}(r)\leq \overline u(r, 0),\;  &0\leq r\leq h_0,
\end{array} \right.
\end{eqnarray*}
then the solution $(u, h)$ of problem $(\ref{a3})$ satisfies $h(t)\leq \overline h(t)$ for $t\in(0, \infty)$
and $u(r, t)\leq \overline u(r, t)$ for $(r, t)\in (0, h(t))\times (0, \infty).$
\end{cor}

\section{The threshold value}

Similarly as the basic reproduction number introduced in epidemiology, in this section, we will give a
threshold value, which plays an important role to determine the spreading or vanishing of the invasive species in a heterogeneous
environment. We will first define a threshold value in a fixed domain $B_R$ (a ball with radius $R$ and center at $0$) and then
 give a threshold value for the free boundary problem (\ref{a3}).

If the environment is heterogeneous and no individual on the boundary, the correspondence model of system (\ref{a3}) with the fixed habitat $B_R$ is
described by
\begin{eqnarray}
\left\{
\begin{array}{lll}
u_t-D\Delta u=b(r)u-d(r)u-\beta(r) u^2,\; &r\in B_R,\;t>0,  \\
u=0,  &r\in \partial B_R ,\;t>0.
\end{array} \right.
\label{B11}
\end{eqnarray}
Define the threshold value $R^{Di}_0(D, B_R)$ for (\ref{B11}) with homogeneous Dirichlet boundary condition by
$$R^{Di}_0(D, B_R):=\ \sup_{\phi\in H^1_0(B_R),\phi\neq 0}\left \{\frac{\int_{B_R} b(r)
\phi^2dr}{\int_{B_R} (D|\triangledown \phi|^2+d(r)\phi^2)dr}\right\}.$$
If $b$ and $d$ are positive constants, then it is easy to check that
$R^{Di}_0(D, B_R)= \frac{b}{\lambda(R)D+d}$, where $\lambda(R)$
is the principle eigenvalue of the following problem
\begin{eqnarray}
\left \{
\begin{array}{lll}
-\Delta\psi=\lambda\psi,\; &r\in B_R,\\
\psi=0,  &r\in \partial B_R.
\end{array} \right.
\end{eqnarray}
It is well known that $\lambda(R)$ is a strictly decreasing continuous function with $R$
and
$$\lim_{R\rightarrow0^{+}}\lambda(R)=+\infty\ \textrm{and}\
\lim_{R\rightarrow+\infty}\lambda(R)=0.$$

With the above definition, we have the following statements:
\begin{thm} The following conclusions about $R_0^{Di}(D,B_R)$ hold:

$(a)$ $R_0^{Di}(D,B_R)$ is a positive and monotone decreasing function of $D$, $R_0^{Di}(D,B_R)\rightarrow 0$ as $D\rightarrow\infty$ and
$R_0^{Di}(D,B_R)\rightarrow \sup_{r\in B_R}\frac{b(r)}{d(r)}$ as $D\rightarrow 0^+$;

$(b)$ $R_0^{Di}(D,B_R)$ is a strictly monotone increasing function of the radius $R$. That is to say,
$R_0^{Di}(D,B_{R_1})<R_0^{Di}(D, B_{R_2})$ if $R_1<R_2$.
Furthermore, $R_0^{Di}(D,B_R)\rightarrow \sup_{r\in[0,\infty)}\frac{b(r)}{d(r)}$ as $R\rightarrow\infty$.

$(c)$ There exists a threshold value $D^*:=D^*(R)\in [0, \infty)$ such that $R_0^{Di}(D,B_R)\geq 1$ for $D\leq D^*$ and $R_0^{Di}(D,B_R)<1$ for $D>D^*$. Moreover,
  $D^*\in (0, \infty)$ if $\{r\in B_R:\, b(r)>d(r)\}$ is nonempty;

 $(d)$ There exists a threshold value $h^*:=h^*(D)\in (0, \infty]$ such that $R_0^{Di}(D,B_R)\geq 1$ for $R\geq h^*$ and $R_0^{Di}(D,B_R)<1$ for $R<h^*$. Moreover,
  $h^*\in (0, \infty)$ if $F^{+}$ is nonempty.
\end{thm}

\bpf
The positivity and monotonicity in $(a)$ can be obtained directly by the definition of $R_0^{Di}(D, B_R)$.
It is easy to see that $R_0^{Di}(D,B_R)\rightarrow 0$ as $D\rightarrow\infty$ from the definition.
The proof for that $R_0^{Di}(D,B_R)\rightarrow \sup_{r\in B_R}\frac{b(r)}{d(r)}$ as $D\rightarrow 0^+$ is similar to that of Lemma 2.3 in \cite{L1}.

The proof of $(b)$ is similar to that of Theorem 3.2(a) in \cite{ZX}. Thus we omit the detail here.
And $(c)$ and $(d)$ can be shown directly from parts $(a)$ and $(b)$, respectively.
According to (3.2) in \cite{CC1}, the threshold value $D^*$ in $(c)$ can be described by
{\small  $$D^*=\ \sup_{\phi\in H^1_0(B_R)} \left\{\frac{\int_{B_R} (b(r)-d(r))
\phi^2dr}{\int_{B_R} |\triangledown \phi|^2dr}: \int_{B_R} (b(r)-d(r))
\phi^2dr>0 \right\}.$$}
\epf

In addition, if $\lambda^*:=\lambda^*(D,B_R)$ is the first eigenvalue of
\begin{eqnarray}
\left\{
\begin{array}{lll}
-D\Delta\varphi=b(r)\varphi-d(r)\varphi+\lambda\varphi,\ &r\in B_R,\\
\varphi=0,&r\in \partial B_R,
\end{array} \right.
\end{eqnarray}
its corresponding eigenfunction $\varphi^*$ can be chosen to be positive on $B_R$.
It is known that $\lambda^*$ is determined by variational characterization:
{\small $$\lambda^*=\inf_{\phi\in H^1_0(B_R)}\left\{\int_{B_R}[D|\nabla\varphi|^2+(d(r)-b(r))\varphi^2]dr:\int_{B_R}\varphi^{2}dr=1\right\}.$$}
With the above defined threshold value $R_0^{Di}(D,B_R)$ and the principle eigenvalue $\lambda^*(D,B_R)$,
we have
\begin{lem}
$1-R_0^{Di}(D,B_R)$ has the same sign as $\lambda^*(D,B_R)$.
\end{lem}

It is easy to see that the habitat $B_{h(t)}$ of the species in the free boundary problem \eqref{a3} is changing with $t$.
Therefore the threshold value $R_0^{Di}(D,B_{h(t)})$ is not a constant but a function with $t$. Let us introduce
the threshold value $R_0^{Fr}(D,t)$ for the free boundary problem \eqref{a3} by
$$R_0^{Fr}(D,t):=R_0^{Di}(D,B_{h(t)})=\ \sup_{\phi\in H^1_0(B_{h(t)}),\phi\neq 0}\left \{\frac{\int_{B_{h(t)}} b(r)
\phi^2dr}{\int_{B_{h(t)}} (D|\triangledown \phi|^2+d(r) \phi^2)dr}\right\}.$$
Then it follows from Lemma 2.2 and Theorem 3.1 that
\begin{thm} The following conclusions about $R_0^{Fr}(D, t)$ hold:

$(1)$ $R_0^{Fr}(D,t)$ is a strictly monotone increasing function of $t$, namely, if $t_1<t_2$, then
we have $R_0^{Fr}(D,t_1)<R_0^{Fr}(D,t_2)$. Moreover, if $h_\infty=\infty$, then
$R_0^{Fr}(D,t)\rightarrow \sup_{r\in[0,\infty)}\frac{b(r)}{d(r)}$ as $t\rightarrow\infty$;

$(2)$ If $\{r\in B_{h(t_0)}:\, b(r)>d(r)\}$ is nonempty or the habitat is favorable
at some $t_0\geq0$, there exists a threshold value $D^*:=D^*(B_{h(t_0)})\in (0, \infty)$ such that $R_0^{Fr}(D,t_0)\geq 1$ if $D\leq D^*$
and $R_0^{Fr}(D,t_0)<1$ if $D>D^*$;

$(3)$ If $F^{+}$ is nonempty, there exists a threshold value $h^*\in (0, \infty)$ such that $R_0^{Fr}(D,t)\geq 1$ for $t$ satisfying $h(t)\geq h^*$ and $R_0^{Fr}(D,t)<1$ for $t$ satisfying $h(t)<h^*$.
\end{thm}

\section{Spreading}

In this section, we will show that the slow invasive species living in a favorable habitat will spreading in the whole domain. For this purpose,
we first prove that if the habitat is limited in the long run, then the species vanishes.
\begin{lem}If $h_\infty<\infty$, then $\lim_{t\rightarrow+\infty}\|u(\cdot,t)\|_{C([0,h(t)])}=0.$
\end{lem}
\bpf
 Assume that
$$\limsup_{t\to+\infty} \ ||u(\cdot, t)||_{C([0, h(t)])}=\delta>0$$
by contradiction.
Then there exists a sequence $\{(r_k,t_k)\}$ in $[0,h(t))\times(0,\infty)$
such that $u( r_k,t_k)\geq \delta /2$ for all $k \in \mathbb{N}$,
and $t_k\to \infty$ as $k\to \infty$.

Now we show that
\begin{eqnarray}\|u\|_{C^{
1+\alpha,(1+\alpha)/2}([0, h(t))\times[0, \infty))}+\|h\|_{C^{1+\alpha/2}([0,\infty))}\leq
M\label{Bg1}
\end{eqnarray}
for any $\alpha\in(0,1)$, where constant $M$ depends on
$\alpha, h_0, \|u_0\|_{C^2([0, h_0])}$ and $h_\infty$.
Straighten the free boundary $r=h(t)$ by the transformation $s=\frac{rh_0}{h(t)}$ to the line $s=h_0$.
Let $v(s,t)=u(r,t)$, then we know that $v(s,t)$ satisfies
\begin{eqnarray}
\left\{
\begin{array}{lll}
v_{t}-\frac{h^2_0}{h^2(t)}D\Delta v-\frac{h'(t)}{h(t)}sv_s=v(b -d-\beta v),\; &0<s<h_0,\; t>0, \\
v_s(0,t)=v(h_0,t)=0,&t>0,\\
v(s,0)=u_0(s)\geq 0,\; & 0\leq s\leq h_0.
\end{array} \right.
\label{Bb1}
\end{eqnarray}
Using Lemma 2.2, there exist constants $M_1,\,M_2$ and $M_3$ such that
$$\|\frac{h^2_0}{h^2(t)}\|_{L^\infty}\leq M_1,\,\|\frac{h'(t)}{h(t)}s\|_{L^\infty}\leq M_2,\,
\|v(b -d-\beta  v)\|_{L^\infty}\leq M_3.$$
Applying the parabolic $L^p$ theory  and the Sobolev embedding
theorem (see \cite{LSU, L}) yields
\begin{eqnarray*}
\|v\|_{C^{1+\alpha,(1+\alpha)/2}([0,
h_0]\times[0,\infty))}\leq  M_4.
\end{eqnarray*}
It is easy to see that (\ref{Bg1}) holds.
Noting that $h(t)$ is bounded and $\|h\|_{C^{1+\alpha/2}([0,\infty))}\leq M$,
we have $h'(t)\rightarrow0\,\textrm{as}\, t\rightarrow+\infty$.
It follows from the free boundary condition that $u_r(h(t),t)$ converges to $0$ as the time $t$ goes to infinity.

On the other hand, since $0\leq r_k<h(t)<h_\infty<\infty$,
then there exists a subsequence $\{r_{k_n}\}$ which converges to $r_0\in [0, h_\infty)$
as $n\rightarrow\infty$. For convenience, we denote $\{r_{k_n}\}$ as $\{r_k\}$.
Thus we have $r_k\rightarrow r_0\in [0, h_\infty)$ as $k\rightarrow\infty$.
Define $U_k(r,t)=u(r,t_k+t)$ for $r\in[0,h(t_k+t)),\;t\in(-t_k,\infty)$.
According to the parabolic regularity, then $\{U_k\}$ has a subsequence $\{U_{k_n}\}$ which
converges to $\widehat{U}$ as $n$ approaches to infinite,
and $\widehat{U}$ satisfies
$$
\widehat{U} _t-D\Delta \widehat{U}=(b -d-\beta\widehat{U})\widehat{U},
\; r\in(0,h_\infty), \ t\in (-\infty, \infty).
$$
Since that $\widehat{U}(r_0,0)\geq\delta/2$, we have $\widehat{U}>0$ in
$[0,h(t))\times(-\infty,\infty)$ by the strong comparison principle.
Let $M=\|b(r)-d(r)-\beta(r)\widehat{U}\|_{L^\infty(0, \infty)}$, then
$\widehat{U} _t-D \Delta \widehat{U}\geq-M\widehat{U}$. Using the Hopf Lemma
at the point $(h_\infty,0)$, we get $\widehat{U}_r(h_\infty, 0 )\leq -\sigma^*$ for some $\sigma^*>0$.

Furthermore, the fact $\|u\|_{C^{1+\alpha,(1+\alpha)/2,}([0, h(t))\times[0, \infty))}\leq M$
implies that $u_r(h(t_k),t_k)=(U_k)_r(h(t_k),0)\to\widehat{U}_r(h_\infty,0)$ as $k\to \infty$.
That is contradiction to that $u_r(h(t),t)\rightarrow 0$ given above.
The proof is complete.
\epf
\begin{lem}
If $R_0^{Fr}(D,0)\left(:=R_0^{Di}(D,B_{h_0})\right)\geq1$, then $h_\infty=\infty$.
\end{lem}
\bpf
$R_0^{Fr}(D, 0)\geq1$ can be considered in two cases. One is $R_0^{Fr}(D,0)>1$, the other is $R_0^{Fr}(D,0)=1$.
Noting that $R_0^{Fr}(D,t)$ is strictly monotone increasing function of $t$ in Theorem 3.3(1).
Therefore, $R_0^{Fr}(D, T_0)>R_0^{Fr}(D,0)=1$ for any positive time $T_0$. Setting $T_0$
be the initial time, then we go back the first case $R_0^{Fr}(D, 0)>1$.
Thus we only show the case of $R_0^{Fr}(D,0):=R_0^{Di}(D,B_{h_0})>1$.
By Lemma 3.2, we have $\lambda^*<0$, where $\lambda^*$ satisfies the following eigenvalue problem
\begin{eqnarray}
\left\{
\begin{array}{lll}
-D\Delta \psi=b(r) \psi-d(r)\psi+\lambda^* \psi,\; &r\in B_{h_0},  \\
\psi=0, &r\in \partial B_{h_0},
\end{array} \right.
\label{B2f}
\end{eqnarray}
and the corresponding eigenfunction $\psi$ is positive and $\|\psi\|_{L^\infty}=1$.

Now we construct the suitable lower solution of problem (\ref{a3}). Setting $B=\|\beta(r)\|_{L^\infty}$ and
defining $$\underline{u}=\sigma\psi\; \textrm{for}\;0<r<h_0,\; t>0,$$
where $\sigma$ is small enough such that $0<\sigma\leq-\frac{\lambda^*}{B}$ and
$\sigma\psi\leq u_0(r)$ in $[0,h_0]$.

 Direct computation shows that
$$\underline{u}_t-D\Delta\underline{u}-(b(r)-d(r)-\beta\underline{u})\underline{u}
=\sigma\psi(\lambda^*+\beta\underline{u})\leq0.$$
Then we obtain
\begin{eqnarray*}
\left\{
\begin{array}{lll}
\underline{u}_t-D\Delta \underline{u}\leq (b(r)-d(r)-\beta(r)\underline {u})\underline {u},\; &0<r<h_0,\ t>0, \\
\underline{u}_{r}(0,t)=0,\ \underline{u}(h_0,t)=0,\,\; & t>0, \\
0=h'_0\leq -\mu \underline{u}_{r}(h_0,t)=-\mu \sigma\psi'(h_0),\  &t>0,\\
\underline {u}(r,0)\leq u_{0}(r),\; &0<r<h_0.
\end{array} \right.
\end{eqnarray*}
Using Lemma 2.4 yields
$$u(r,t)\geq\underline{u}(r,t) \;\textrm{for}\; (r,t)\in [0,h_0]\times[0,\infty).$$
It follows that $\liminf_{t\rightarrow+\infty}u(0,t)\geq\sigma\psi(0)>0$.
By Lemma 4.1, we know that $h_\infty=\infty.$
\epf

Assume that $F^+$ is nonempty, it follows from Theorem 3.1(d) that there exists $h^*\in (0, \infty]$ satisfying $R_0^{Di}(D, B_{h^*})=1$. If $h_0$
satisfies $h_0\geq h^*$, we have $R_0^{Di}(D,B_{h_0})\geq R_0^{Di}(D,B_{h^*})=1$ from Theorem 3.1 (d).
The following result follows from Lemma 4.2 directly.
\begin{cor}If $h_0\geq h^*$, then $h_\infty=\infty$.
\end{cor}

In order to study the long time behavior of the spreading species, we consider the fixed boundary problem
\begin{eqnarray}
\left\{
\begin{array}{lll}
u_{t}-D\Delta u=(b(r)-d(r)-\beta(r) u)u,\; & 0<r<R,\,t>0,  \\
u_r(0,t)=u(R,t)=0,&t>0,\\
u(r,0)=u_0(r)\geq0, \ & 0<r<R,
\end{array} \right.
\label{111c}
\end{eqnarray}
where $R$ is a positive number and $u_0(r)$ is not always equal to $0$.
And the related stationary problem is
\begin{eqnarray}
\left\{
\begin{array}{lll}
-D\Delta u=(b(r)-d(r)-\beta(r) u)u,\; & 0<r<R,  \\
u_r(0)=u(R)=0,\\
u_0(r)\geq0, \ & 0<r<R.
\end{array} \right.
\label{11c}
\end{eqnarray}
For problems $(\ref{111c})$ and $(\ref{11c})$, the following results are known.
\begin{prop}
$(a)$ If $R_0^{Di}(D,B_R)>1$, then the elliptic problem in a fixed domain $(\ref{11c})$ has only one
positive solution $\widetilde{u}_R(r)$.

$(b)$ If $R_0^{Di}(D,B_R)>1$. Then any positive solution of problem $(\ref{111c})$
converges to $\widetilde{u}_R(r)$ as the time approaches to infinity.
\end{prop}
The proof of $(a)$ is exactly the same as Theorem 2.1 in \cite{Dong}
or \cite{HH} and
 $(b)$ can be seen in \cite{DL, HH}. So we omit the detail here.
\begin{prop}
Assume $(H)$ holds. Then

$(1)$ the positive solution $\widetilde{u}(r)$ of problem $(\ref{b121212})$ is unique;

$(2)$ $\widetilde{u}_R(r)\rightarrow \widetilde{u}(r)$ as $R \rightarrow +\infty$,
where $\widetilde{u}_R(r)$ is the unique solution to problem $(\ref{11c})$.
\end{prop}
For the proof, $(1)$ follows from Theorem 1.2 in \cite{Dong} and $(2)$ follows from
Proposition 3.4 in \cite{Dong}.

\begin{lem} Suppose that $(H)$ holds. If $h_\infty=\infty$. Then
$\lim_{t\rightarrow+\infty}u(r,t)=\widetilde{u}(r)$ locally uniformly in $[0,+\infty)$,
where $\widetilde{u}(r)$ is the unique positive solution of problem $(\ref{b121212})$.
\end{lem}
\bpf Since $\lim_{r\rightarrow\infty}(b(r)-d(r))=\alpha$, there exists a large positive constant
$r^*$ such that $b(r)-d(r)>0$ and $R_0^{Di}(D,B_r)>1$ for any $r>r^*$ by Theorem 3.1(b). And we note that there exists a $t^*>0$ such that $r^*=h(t^*)$.

Next we choose increasing sequence $\{R_n\}$ such that $\lim_{n\rightarrow+\infty}R_n=+\infty$,
$R_1>r^*$ and $R_0^{Di}(D,B_{R_n})>1$ for any $n$. Therefore, for each fixed $n$,
we can find $t_{n}>0$ such that $h(t_{n})=R_{n}$. Let $u_{n}(r,t)$ be the positive
solution of the following problem
\begin{eqnarray*}
\left\{
\begin{array}{lll}
w_{t}-D\Delta w=(b(r)-d(r)-\beta(r) w)w,\; & 0<r<R_{n},\,t>t_{n},  \\
w_r(0,t)=w(R_{n},t)=0,&t>t_{n},\\
w(r,t_{n})=u(r,t_{n})\geq0, \ & 0<r<R_{n}.
\end{array} \right.
\end{eqnarray*}
Using the comparison principle, we obtain
$$u_{n}(r,t)\leq u(r,t)\;\textrm{for}\;r\in[0,R_{n}],\, t\geq t_{n}.$$
By virtue of Proposition 4.4(b),
$$\widetilde{u}_{R_{n}}(r)\leq \underline{\lim}_{t\rightarrow+\infty}u(r,t)\; \textrm{uniformly in} \;[0,R_{n}].$$
Thus, for any $L>0$, there exists $n_0$ such that $R_n\geq L$ for $n\geq n_0$.
We then have
$$\widetilde{u}_{R_{n}}(r)\leq \underline{\lim}_{t\rightarrow+\infty}u(r,t)\; \textrm{uniformly in} \;[0,L]$$
for $n\geq n_0$. Now let $n\rightarrow+\infty$, Proposition 4.5(2) implies that
\begin{equation}\label{upper}
\widetilde{u}(r)\leq \underline{\lim}_{t\rightarrow+\infty}u(r,t)\; \textrm{locally uniformly in} \;[0,+\infty).
\end{equation}

Next, motivated by the squeezing technique in  \cite{DM}, we claim
$$\overline{\lim}_{t\rightarrow+\infty}u(r,t)\leq\widetilde{u}(r) \; \textrm{locally uniformly in} \;[0,+\infty)$$
by considering the boundary blow-up problem
\begin{eqnarray}
\left\{
\begin{array}{lll}
-D\Delta v^*=(b(r+r^*)-d(r+r^*))v^*-\beta(r+r^*){v^*}^2,\ &0<r<R,\\
v^*|_{\partial B_R}=\infty,
\end{array} \right.
\label{1}
\end{eqnarray}
where $v_R$ is the solution of problem (\ref{1}). It follows from the comparison principle Lemma 2.1 in \cite{DM} that
$v_R$ decreases to $\widetilde{u}$ as $R$ approaches to $\infty$.

Denote $v^*_{n}(r)=v^*(r+r^*)$ be the positive solution of the problem (\ref{1}) with $R$ replaced by $R_{n}$.
And $v^*_{n}(r)=\lim_{i\rightarrow\infty}v^*_{n, i}(r)$ is the unique nonnegative solution of the following problem
\begin{eqnarray*}
\left\{
\begin{array}{lll}
-D\Delta v=(b(r+r^*)-d(r+r^*))v-\beta(r+r^*)v^2,\ &0<r<R_{n},\\
v|_{\partial B_{R_{n}}}=i,
\end{array} \right.
\end{eqnarray*}
where $i$ is the positive integer.

Moreover, we know that $v^*_{n,i}(r)$ is the stationary state of the parabolic problem
\begin{eqnarray}
\left\{
\begin{array}{lll}
\nu_{t}-D\Delta \nu=(b(r+r^*)
-d(r+r^*)-\beta(r+r^*) \nu)\nu,\; & 0<r<R_{n},\,t>0,  \\
\nu(r+r^*,t^*)=k_{n}v^*_{n+1}(r),&0<r<R_{n},\\
\nu(r+r^*,t+t^*)=i, \ & r\in\partial{B_{R_{n}}},\,t>0,
\end{array} \right.
\label{equa}
\end{eqnarray}
where $k_{n}$ is a positive constant and satisfies $k_{n}v^*_{n+1}|_{\partial B_{R_{n}}}\geq i$.
Denote $v_{n}(r,t)=\nu (r+r^*, t+t^*)$ be the solution of problem (\ref{equa}).
Similarly as in \cite{Sa}, we can deduce that $v_{n}(r,t)\rightarrow v_{n,i}^*(r)$ locally uniformly as $t \rightarrow \infty$.
Then it follows from Theorem 4.1 in \cite{DM} that
$$u(r+r^*,t+t^*)\leq v_{n}(r,t)\;\mbox{for}\; 0<r<R_{n}\;\mbox{and}\;t>0,$$
which implies
$$\overline{\lim}_{t\rightarrow\infty}u(r+r^*,t+t^*)\leq v^*_{n,i}(r)
\;\mbox{locally uniformly in}\; [0,R_{n}].$$
Therefore, we get
$$\overline{\lim}_{t\rightarrow+\infty}u(r,t)\leq\widetilde{u}(r) \; \textrm{locally uniformly in} \;[0,+\infty),$$
which together with (\ref{upper}) gives the desired result.
\epf

If $B_{h(t_0)}$ is a favorable habitat and the diffusion is slow. Applying Theorem 3.3(2), we have $R_0^{Fr}(D,t)\geq 1$ for $t\geq t_0$.
Similarly, if $B_{h(t_0)}$ is a favorable habitat and the occupying habitat is large, we then have $R_0^{Fr}(D,t)\geq 1$.
Therefore, it follows from Lemmas 4.2 and 4.6 that the species is spreading.
\begin{thm}
Suppose that $(H)$ holds and $B_{h(t_0)}$ is a favorable habitat for some $t_0\geq0$. If the diffusion is slow $(D\leq D^*(B_{h(t_0)}))$ or the occupying habitat
is large $(h(t_0)\geq h^*$ with $R_0^{Di}(D,B_{h^*})=1)$, then $h_\infty=\infty$
and $\lim_{t\rightarrow+\infty}u(r,t)=\widetilde{u}(r)$ locally uniformly in $[0,+\infty)$,
where $\widetilde{u}(r)$ is the unique positive solution of problem $(\ref{b121212})$.
\end{thm}

\section{The spreading-vanishing dichotomy}
This section is contributed to the spreading-vanishing dichotomy which is
governed by the initial value $u_0(r)=\delta\varphi(r)$, where $\delta$ is a positive constant.
And sufficient conditions are given in this part to determine the spreading
 or vanishing of the invasive species .

\begin{lem}
Assume that $F^+$ is nonempty and $R_0^{Fr}(D,0)<1$.
Then $h_\infty=\infty$ if $||u_0(r)||_{C([0,h_0])}$ is sufficiently large.
\end{lem}
\bpf
Noting that $\lim_{T\rightarrow\infty}R_0^{Di}(D,B_{\sqrt{T}})
=\sup_{r\in[0,\infty)}\frac{b(r)}{d(r)}>1$ by Theorem 3.1(b). Then there exists a $T^*$ such that
$R_0^{Di}(D, B_{\sqrt{T^*}})>1$.

Next we construct a suitable lower solution $(\underline{u}(r,t),\underline{h}(t))$ to problem (\ref{a3}),
where $\underline{h}(t)$ satisfies
$$R_0^{Fr}(D,T^*)=R_0^{Di}(D,B_{h(T^*)})\geq R_0^{Di}(D,B_{\underline{h}(T^*)})\geq R_0^{Di}(D,B_{\sqrt{T^*}})>1.$$

Let $\lambda$ be the eigenvalue of
\begin{eqnarray}
\left\{
\begin{array}{lll}
-D\Delta \varphi-\frac{1}{2}\varphi'=\lambda\varphi,\; & 0<r<1,  \\
\varphi'(0)=\varphi(1)=0,\\
\end{array} \right.
\end{eqnarray}
the corresponding function $\varphi>0$ and $\varphi'>0$ in $[0,1)$, and $\|\varphi\|_{{L^\infty}[0,1)}=1$.

Defining
$$\underline{h}(t)=\sqrt{t+\delta},\;t\geq0,$$
$$\underline{u}(r,t)=\frac{M}{(t+\delta)^k}\varphi(\frac{r}{\sqrt{t+\delta}}),
\;0\leq r\leq\sqrt{t+\delta},\;t\geq0,$$
where $\delta,\;k,\;M$ are positive constants to be chosen later, we will show that $(\underline u,\underline{h})$ is the lower solution to problem (\ref{a3})

It follows from Lemma 2.2 that there exists a positive constant $C_1$
such that for $r\in[0,h(t)]$ and $t\in[0,T^*]$, we have
$0\leq u(r,t)\leq C_1$. Noting that $b_1\leq b(r),\;d(r),\;\beta(r)\leq b_2$ for $r\in[0,+\infty)$,
then there exists a positive constant $L$ such that $f(r,u)\geq -Lu$,
where $f(r,u):=u[b(r)-d(r)-\beta(r)u]$.
Choosing $0<\delta\leq \min \{1,h^2_0\}$, $k>\lambda+L(T^*+1)$ and
$-2\mu M\varphi'(1)>(T^*+1)^k$, we obtain
\begin{eqnarray*}
& &\underline{u}_t-D\Delta\underline{u}-(b(r)-d(r)-\beta(r)\underline{u})\underline{u}\\
& &\leq -\frac{M}{(t+\delta)^{k+1}}[k\varphi+\frac{r}{2\sqrt{t+\delta}}\varphi'+D\triangle\varphi-
L(t+\delta)\varphi]\\
& &\leq-\frac{M}{(t+\delta)^{k+1}}[D\triangle\varphi+\frac{1}{2}\varphi'+k\varphi-L(t+\delta)\varphi]\\
& &\leq-\frac{M}{(t+\delta)^{k+1}}(D\triangle\varphi+\frac{1}{2}\varphi'+\lambda\varphi)=0
\end{eqnarray*}
for  $0<r<\underline{h}(t)$ and $0<t\leq T^*$.
$$\underline{h}'(t)+\mu \underline{u}_r(\underline{h}(t),t)
=\frac{1}{2\sqrt{t+\delta}}
+\frac{\mu M}{(t+\delta)^{k+\frac{1}{2}}}\varphi'(1)<0\;\textrm{for}\;0<t\leq T^*.$$
Moreover, if $\underline{u}(r,0)=\frac{M}{\delta^k}\varphi(\frac{r}{\sqrt{\delta}})<u_{0}(r)$ in $[0,\sqrt{\delta}]$, then we have
\begin{eqnarray*}
\left\{
\begin{array}{lll}
\underline{u}_t-D\Delta \underline{u}\leq (b(r)-d(r)-\beta(r)\underline {u})\underline {u},\; &0<r<\underline{h}(t),\ 0<t\leq T^*,\\
\underline{u}_{r}(0,t)=0,\ \underline{u}(\underline{h}(t),t)=0,\,\; &0<t\leq T^*, \\
\underline{h}'(t)< -\mu \underline{u}_{r}(\underline{h},t),\  &0<t\leq T^*,\\
\underline{u}(r,0)<u_{0}(r),\; &0\leq r\leq\sqrt{\delta}.
\end{array} \right.
\end{eqnarray*}
Noting that $\underline{h}(0)=\sqrt{\delta}\leq h_0$, we can use Lemma 2.4 to conclude that
$h(t)\geq\underline{h}(t)$ in $[0,T^*]$. Specially, we obtain
$h(T^*)\geq\underline{h}(T^*)=\sqrt{T^*+\delta}\geq\sqrt{T^*}$.
Therefore, $R_0^{Fr}(T^*)=R_0^{Di}(D, B_{h(T^*)})\geq R_0^{Di}(D,B_{\sqrt{T^*}})>1$, which together with Lemma 4.2 gives that $h_\infty=\infty$.
\epf

Corollary 4.3 shows that if $h_0\geq h^*$ with $R_0^{Di}(D,B_{h^*})=1$, spreading must happen.
For the case $h_0<h^*$,  the following result is directly from Theorem 3.1(d) and Lemma 5.1.
\begin{cor} Suppose that the favorable set $(F^+)$ is nonempty and $h_0<h^*$. Then $h_\infty=\infty$
if $||u_0(r)||_{C([0, h_0])}$ is sufficiently large.
\end{cor}

\begin{lem}Assume $R_0^{Fr}(D,0)\left(:=R_0^{Di}(D,B_{h_0})\right)<1$,
then $$h_\infty<\infty\;\mbox{and}\; \lim_{t\rightarrow+\infty}\|u(\cdot,t)\|_{C([0,h(t)])}=0,$$
if $\|u_0\|_{C([0,h_0])}$ is sufficiently small.
\end{lem}
\bpf
Since that $R_0^{Di}(D, B_{h_0})<1$, we have $\lambda^*>0$ which satisfies
\begin{eqnarray*}
\left\{
\begin{array}{lll}
-D\Delta\psi=b(r)\psi-d(r)\psi+\lambda\psi,\ &0<r<h_0,\\
\psi_r(0)=\psi(r)=0,&r=h_0,
\end{array} \right.
\end{eqnarray*}
where the eigenfunction $\psi$ is positive in $[0,h_0)$.

Inspired by \cite{DL}, we define
$$\sigma(t)=h_0(1+\sigma-\frac{\sigma}{2}e^{-\sigma t})\;\textrm{for}\;t\geq0,$$
$$v(r,t)=\varepsilon e^{-\sigma t}\psi(\frac{rh_0}{\sigma(t)})\;\textrm{for}\;0\leq r\leq\sigma(t),\;t\geq0,$$
where $\sigma$ and $\varepsilon$ are small positive constants which will be determined below.

For $0<r<\sigma(t)$ and $t>0$, direct calculation gives
\begin{eqnarray*}
& &v_t-D\Delta v-(b-d)v\\
& &=-\varepsilon \sigma e^{-\sigma t}\psi-\frac{rh_0\sigma'(t)}{\sigma^2(t)}\varepsilon e^{-\sigma t}\psi'
-\frac{Dh_0^2}{\sigma^2(t)}\varepsilon e^{-\sigma t}\Delta\psi-(b-d)\varepsilon e^{-\sigma t}\psi\\
& &\geq-\sigma v+\frac{h^2_0}{\sigma^2(t)}(b -d+\lambda^*)v-(b-d)v\\
& &=v[-\sigma+(\frac{h^2_0}{\sigma^2(t)}-1)(b-d)+\frac{\lambda^*h^2_0}{\sigma^2(t)}]\\
& &\geq v[-\sigma+(\frac{1}{(1+\sigma)^2}-1)b+\frac{\lambda^*}{(1+\sigma)^2}]\geq 0
\end{eqnarray*}
if $\sigma$ is sufficiently small such that
$\sigma(1+\sigma)^2+[(1+\sigma)^2-1]b_2\leq\lambda^*$.

Moreover, $\sigma'(t)\geq-\mu v_r(\sigma(t),t)$ if we choose $0<\varepsilon\leq-\frac{\sigma^2h_0}{2\mu\psi'(h_0)}$.
And $u_0$ is small such that $\|u_0\|_{C([0,h_0])}\leq\varepsilon\psi(\frac{h_0}{1+\frac{\sigma}{2}})$, then we have
\begin{eqnarray*}
\left\{
\begin{array}{lll}
v_t-D\Delta v\geq(b(r)-d(r))v,\; &0<r<\sigma(t),\ t>0,\\
v_{r}(0,t)=0,\ v(\sigma(t),t)=0,\,\; & t>0, \\
\sigma'(t)>-\mu v_{r}(\sigma(t),t),\  & t>0,\\
u_{0}(r)\leq\varepsilon\psi(\frac{h_0}{1
+\frac{\sigma}{2}})\leq v(r,0),\; &0\leq r\leq h_0,\\
\sigma(0)=(1+\frac{\sigma}{2})h_0>h_0.
\end{array} \right.
\end{eqnarray*}
Applying Corollary 2.5, we get
$h(t)\leq\sigma(t)$ and $u(r,t)\leq v(r,t)$ for $0\leq r\leq h(t)$, $t>0$.
Furthermore, we have $h_\infty\leq \lim_{t\rightarrow\infty}\sigma(t)
=h_0(1+\sigma)<\infty$. Using Lemma 4.1, it follows that $\lim_{t\rightarrow+\infty}
\|u(\cdot,t)\|_{C([0,h(t)])}=0$.
\epf

Assume that $h_0<h^*$ with $R_0^{Di}(D,B_{h^*})=1$, then
$R_0^{Fr}(D,0)\left(:=R_0^{Di}(D,B_{h_0})\right)<R_0^{Di}(D,B_{h^*})=1$ by Theorem 3.1(b).
The following conclusion comes from Lemma 5.3 directly.
\begin{cor}
Suppose that $h_0<h^*$, then $h_\infty<\infty$
and $\lim_{t\rightarrow+\infty}\|u(\cdot,t)\|_{C([0,h(t)])}=0$ if $\|u_0\|_{C([0,h_0])}$
is sufficiently small.
\end{cor}

If $h_0$ is fixed, the spreading or vanishing of a invasive species depends on the diffusion rate $D$ and the initial number $u_0(r)$ of the species.
\begin{thm} Suppose $(H)$ holds. If $\{r\in B_{h_0}:\, b(r)>d(r)\}$ is nonempty, then there exists $D^*\in (0,\infty)$ such that

$(i)$ the spreading happens for any initial value if $D\leq D^*$;

$(ii)$ when $D>D^*$, then one of the following results hold: Vanishing occurs if $\|u_0\|_{C([0,h_0])}$is sufficiently small;
Spreading occurs if $\|u_0\|_{C([0,h_0])}$ is sufficiently large.
\end{thm}
\bpf
According to Theorem 3.3(2), we know that there exists a constant $D^*\in (0,\infty)$
such that $R_0^{Fr}(D,0)\geq 1$ if $D\leq D^*$.
In the case of $R_0^{Fr}(D,0)\geq 1$, using Lemmas 4.2 and
 4.6 we can conclude that the invasive species is spreading.
Thus $(i)$ can be obtained.

On the other hand, if $D>D^*$, then we have $R_0^{Fr}(D,0)<1$.
If $(H)$ holds, the favorable set $F^+\neq\emptyset$.
So $(ii)$ follows from Lemmas 5.1 and 5.3 directly. The proof is complete.
\epf

Similarly, if $D$ is fixed, the spreading or vanishing of a invasive species depends on the radius $h_0$ of the initial occupying
habitat $B_{h_0}$ and the initial number $u_0(r)$ of the species.
\begin{thm} Suppose $(H)$ holds. There exists $h^*\in (0,\infty)$ such that

$(i)$ the spreading happens for any initial value if $h_0\geq h^*$;

$(ii)$ when $h_0<h^*$, then one of the following results hold: Vanishing occurs if $\|u_0\|_{C[0,h_0]}$ is sufficiently small;
Spreading occurs if $\|u_0\|_{C([0,h_0])}$ is sufficiently large.
\end{thm}

If $D$ and $h_0$ are fixed, the initial number $u_0(r)$ governs the spreading and vanishing of the invasive species.
\begin{thm}
Assume that $(H)$ holds.
Let $(u(r,t);h(t))$ be the solution of problem $(1.4)$ with $u_0(r):=
\delta\varphi(r)$ for some constant $\delta>0$.
Then there exists $\delta_0:=\delta_0(h_0;D)\in [0,+\infty)$
such that spreading occurs for $\delta>\delta_0$ and vanishing happens for $0<\delta\leq\delta_0$.
\end{thm}
\bpf
Since $\lim_{r\rightarrow\infty}(b(r)-d(r))=\alpha$, then there exists
$ h^*(D)\in(0,\infty)$ such that $R_0^{Di}(D;B_{h^*})=1$ by Theorem 3.1 (d).
Applying Corollary 4.3 and Lemma 4.6, we know that spreading happens
if $h_0\geq h^*(D)$. Then we can choose $\delta_0=0$ in this case.

Next, we consider the other case $h_0<h^*$. Thus we have $R_0^{Di}(D,B_{h_0})<1$ from
Theorem 3.1 (d). Define
$$\delta_0:=\sup\left\{\delta:\;h_\infty<\infty\;\mbox{for}\;\delta^0\in(0,\delta]\right\}.$$
It follows from Corollary 5.4 that vanishing occurs for small $\delta>0$.
Hence $\delta_0 \in (0, +\infty]$. Moreover, Lemma 5.1 illustrates that spreading happens for
large $\delta$. Therefore $\delta_0 \in (0, +\infty)$. According to the comparison principle,
we deduce that spreading occurs when $\delta>\delta_0$; while if $0<\delta<\delta_0$, the vanishing happens.

Finally, we show that vanishing occurs when $\delta=\delta_0$.
Otherwise, spreading must happen and we have $h_\infty=\infty$ at $\delta=\delta_0$.
It follows from Theorem 3.3(1) that there exists $T_0>0$ such that
$R_0^{Fr}(T_0)=R_0^{Di}(D,B_{h(T_0)})>1$. To stress the continuous dependence of $(u,h)$
on the initial value $u_0(r)=\delta\varphi(r)$, let $(u_\varepsilon;h_\varepsilon)$
denote the solution of problem (1.4) with $u_0(r)=\delta_0\varphi(r)$ replaced by
$u^\varepsilon_0(r)=(\delta_0-\varepsilon)\varphi(r)$. Obviously, we have
$R_0^{Di}(D,B_{h^\varepsilon(T_0)})>1$. But applying Lemma 4.2 yields that
$(u_\varepsilon; h_\varepsilon)$ spreads to the whole habitat.
But this is contradict to the definition of $\delta_0$.
\epf

\section{Spreading speed}

In this section, we always suppose $(H)$ holds. In the spreading case, we will give the asymptotic spreading speed of the free boundary $x=h(t)$.
First, it follows from Section 4 of \cite{DL} that $k_0$ increases with $\mu, a$ and decreases with $b$,
where $k_0\in (0,2\sqrt{aD})$ is the unique positive solution of the following problem
\begin{eqnarray}
\left\{
\begin{array}{lll}
-D U''+kU'=aU-bU^2,\ &r>0,\\
U(0)=0.
\end{array} \right.
\end{eqnarray}
It follows from a well known conclusion that the problem of problem (6.1) has a unique positive solution
$U=U_{a, b, k}$. And $U(r)$ converges to $\frac{a}{b}$ as $r$ approaches to infinity. Moreover, $U'_{a, b, k}(r)$ is
a monotone increasing function about $r(\in[0,+\infty))$, $U'_{a, b, k}(0)$ and $U_{a, b, k}(r)$ are strictly
decreasing functions with $k$ for $r>0$, respectively.

\begin{thm}
If $h_\infty=\infty$, then $\lim_{t\rightarrow+\infty}\frac{h(t)}{t}=k_0$.
\end{thm}
\bpf
Now we show that $u(r,t)$ is locally uniformly bounded for $[0,+\infty)$.
Since the condition (H) holds,
by Theorem 3.1 in \cite{DM}, the unique positive solution $\widetilde{u}(r)$ of problem (1.6) satisfies that
$\lim_{r\rightarrow+\infty}\widetilde{u}(r)=\frac{\alpha}{\beta}$.
For any $\varepsilon>0$, there exists $R:=R(\varepsilon)$ such that for any $r\geq R$,
$$\alpha-2\varepsilon\leq b(r)-d(r)\leq\alpha+2\varepsilon,\;\;
\beta-2\varepsilon\leq\beta(r)\leq\beta+2\varepsilon$$ and
$\frac{\alpha-\varepsilon}{\beta+\varepsilon}\leq\widetilde{u}(r)
\leq\frac{\alpha+\varepsilon}{\beta+\varepsilon}$. As $h_\infty=+\infty$ and
$\lim_{t\rightarrow+\infty}u(r,t)=\widetilde{u}(r)$, then there exists $T:=T(R)>0$ such that
for any $t\geq0$,
$$h(T)>3R \;\;\mbox{and}\;\;u(2R,t+T)<\frac{\alpha+2\varepsilon}{\beta-2\varepsilon}.$$

Let $$\widehat{u}(r,t)=u(r+2R,t+T)\;\;\mbox{and}\;\;\widehat{h}(t)=h(t+T)-2R.$$
Direct calculations yield
\begin{eqnarray}
\left\{
\begin{array}{lll}
\widehat{u}_{t}-D \Delta(\widehat{u}_{rr}+\frac{n-1}{r+2R}\widehat{u}_{r})=&\\
\quad \widehat{u}(b(r+2R) -d(r+2R)-\beta (r+2R)\widehat{u}),\; &0<r<\widehat{h}(t),\; t>0,   \\
\widehat{u}(0,t)=u(2R,t+T),\;\widehat{u}(\widehat{h}(t),t)=0,\;&t>0,\\
\widehat{h}'(t)=-\mu \widehat{u}_{r}(\widehat{h}(t),t),\;   & t>0, \\
\widehat{u}(r,0)=u_{0}(r+2R,T), & 0\leq r\leq \widehat{h}(0).
\end{array} \right.
\label{6012}
\end{eqnarray}
Since $b(r+2R) -d(r+2R)\leq\alpha+2\varepsilon$ and $\beta (r+2R)\geq\beta-2\varepsilon$,
it follows from the comparison principle that
$$\widehat{u}(r,t)\leq w(t)\;\mbox{for}\; 0<r<\widehat{h}(t)\;\mbox{and}\;t>0,$$
where $w(t)$ is the solution of
\begin{eqnarray*}
\left\{
\begin{array}{lll}
\frac{dw}{dt}=w[(\alpha+2\varepsilon)-(\beta-2\varepsilon)w],\ &t>0,\\
w(0)=\max\left\{\frac{\alpha+2\varepsilon}{\beta-2\varepsilon},\;\|\widehat{u}(0,\cdot)\|_\infty\right\}.
\end{array} \right.
\end{eqnarray*}
Noting that $\lim_{t\rightarrow\infty}w(t)=\frac{\alpha+2\varepsilon}{\beta-2\varepsilon}$,
for any $0\leq r\leq\widehat{h}(t)$ and $t\geq \widehat{T}$, there exists
$\widehat{T}:=\widehat{T}_\varepsilon$ such that
$$\widehat{u}(r,t)\leq w(t)\leq \frac{\alpha+2\varepsilon}{\beta-2\varepsilon}(1-\varepsilon)^{-1}.$$
That is to say
$$u(r+2R,t+T)\leq\frac{\alpha+2\varepsilon}{\beta-2\varepsilon}(1-\varepsilon)^{-1}
\;\mbox{for}\; 0\leq r\leq\widehat{h}(t)\;\mbox{and}\;t\geq \widehat{T}.$$

On the other hand, we will claim that $u$ has a subbound.
Replacing $b(r+2R) -d(r+2R)$ and $\beta (r+2R)$ by $\alpha-2\varepsilon$
and $\beta+2\varepsilon$, respectively, then problem (6.2) becomes
\begin{eqnarray*}
\left\{
\begin{array}{lll}
v_{t}-D \Delta(v_{rr}+\frac{n-1}{r+2R}v_{r})=v[(\alpha-2\varepsilon)-(\beta+2\varepsilon)v],\; &0<r<\widehat{h}(t),\; t>0,   \\
v(0,t)=\widehat{u}(0,t),\;v(\widehat{h}(t),t)=0,\;&t>0,\\
v(r,0)=\widehat{u}(r,0), & 0\leq r\leq \widehat{h}(0),
\end{array} \right.
\end{eqnarray*}
It follows from Lemma 3.4 in \cite{DG} that
$$\lim_{t\rightarrow+\infty}v(r,t)\geq\frac{\alpha-2\varepsilon}{\beta+2\varepsilon}\ \mbox{locally uniformly for}\; r\in[0,\infty).$$
Using the comparison principle gives
$$\widehat{u}(r,t)\geq v(r,t)\; \mbox{for}\;r\in[0,\widehat{h}(t)]\;\mbox{and}\;t>0.$$
Thus we conclude that
$$\underline{\lim}_{t\rightarrow\infty}\widehat{u}(r,t)\geq\frac{\alpha-2\varepsilon}{\beta+2\varepsilon}
\ \mbox{locally uniformly for}\;\; r\in[0,\infty).$$

Next, we construct the suitable lower and upper solution of the free boundary problem (1.4) to show that
$\lim_{t\rightarrow+\infty}\frac{h(t)}{t}=k_0$.
The rest part of the proof is similar as that of Theorem 3.6 in \cite{DG} with some obvious modifications. Define
$$\underline{h}(t)=(1-\varepsilon)^2k_0(\mu,\alpha-2\varepsilon,\beta+2\varepsilon)t+\widehat{h}(0),\,t\geq0,$$
$$\underline{u}(r,t)=(1-\varepsilon)^2U_
{\alpha-2\varepsilon,\beta+2\varepsilon,k_0(\mu,\alpha-2\varepsilon,\beta+2\varepsilon)}(\underline{h}(t)-r),
0\leq r\leq\underline{h}(t),\,t>0.$$
Noting that $U_{\alpha+2\varepsilon,\beta-2\varepsilon,k_0(\mu,\alpha+2\varepsilon,\beta-2\varepsilon)}$
converges to $\frac{\alpha+2\varepsilon}{\beta-2\varepsilon}$ as $r$
approaches to infinity. Then there
exists $R_0(\varepsilon)>2R$ such that
$U_{\alpha+2\varepsilon,\beta-2\varepsilon,k_0(\mu,\alpha+2\varepsilon,\beta-2\varepsilon)}
>\frac{\alpha+2\varepsilon}{\beta-2\varepsilon}(1-\varepsilon)$ for $r\geq R_0$.
Let
$$\overline{h}=(1-\varepsilon)^{-2}k_0(\mu,\alpha+2\varepsilon,\beta-2\varepsilon)t+R_0+\widehat{h}(\widehat{T}),\,t\geq0,$$
$$\overline{u}(r,t)=(1-\varepsilon)^{-2}
U_{\alpha+2\varepsilon,\beta-2\varepsilon,k_0(\mu,\alpha+2\varepsilon,\beta-2\varepsilon)}(\overline{h}(t)-r),
0\leq r\leq\overline{h}(t),\,t>0.$$
It follows from the proof of Theorem 3.6 in \cite{DG} and Lemma 2.4, we find that
$(\underline{u},\underline{h})$ and $(\overline{u},\overline{h})$ are the lower and upper solutions of problem
(1.4).

Hence we have
$$\underline{\lim}_{t\rightarrow\infty}\frac{h(t)}{t}\geq\lim_{t\rightarrow\infty}\frac{h(t)}{t}=
(1-\varepsilon)^2k_0(\mu,\alpha-2\varepsilon,\beta+2\varepsilon),$$
$$\overline{\lim}_{t\rightarrow\infty}\frac{h(t)}{t}\leq\lim_{t\rightarrow\infty}\frac{\overline{h}(t)}{t}=
(1-\varepsilon)^{-2}k_0(\mu,\alpha+2\varepsilon,\beta-2\varepsilon).$$
Let $\varepsilon\rightarrow0$, then we get $\lim_{t\rightarrow+\infty}\frac{h(t)}{t}=k_0$.
\epf

\section{Discussion}

In this paper, a logistic reaction-diffusion equation is investigated
in heterogeneous environments with a free boundary describing the front of an invasive species.
Based on the sign of $\int_{B_{h(t)}}(b(r)-d(r))dr$, the
habitat of the invasive species are divided into favorable habitat and unfavorable
habitat. We study the long time behavior of the solution and discuss how the spatial heterogeneity affects
the moving patterns of the invasive species.
Sufficient conditions are given to ensure that the spreading and vanishing happen. Furthermore, when the species survives and establishes
itself successfully in the new environment,
we estimate the asymptotic spreading speed, which is smaller than the minimal speed of the corresponding
traveling wave problem.

For simplicity, we always assume that $(H)$ holds. In a favorable habitat $B_{h(t_0)}$ for some $t_0\geq 0$, if the diffusion
is slow or the the occupying habitat is large, we have $h_\infty=\infty$ and $\lim_{t\rightarrow+\infty}u(r,t)=\widetilde{u}(r)$
locally uniformly in $[0,+\infty)$, where $\widetilde{u}(r)$ is the unique positive solution of problem $(\ref{b121212})$ (Theorem 4.7).
That is to say, if the average birth rate of a species is greater than the average death rate, the invasive species with
slow diffusion or large habitat occupation will survive in the new environment. In a biological view, the species will survive
easily in a favorable habitat.

An unfavorable habitat is bad for the species with small number at the beginning (Lemma 5.3), the rare (or endangered) species
in an unfavorable habitat will become extinct in the future.
However, even the habitat is unfavorable, if the initial occupying area $B_{h_0}$ is beyond a
critical size, namely, $h_0\geq h^*$,
then regardless of the initial population size $u_0(r)$, spreading always happens (Corollary 4.3).
And if $h_0<h^*$, spreading is also possible for big initial population size $u_0(r)$ provided that the favorable set $F^+$
is nonempty (Corollary 5.2). Those results tell us that we can also choose a proper initial
habitat or keep sufficient number to preserve the endangered species.
Theorems 5.5, 5.6 and 5.7 imply that slow diffusion, large occupying habitat and big initial population number are benefit
for the species to survive in the new environment. The invasive species with slow diffusion will spread in the total habitat,
while if the ability of migratory is big, the vanishing or spreading of the species in the new environment is determined by the initial number.
A threshold value about dispersal is given in Theorem 5.5. Theorem 5.6
also shows a similar result about the initial occupying habitat.
However, initial value also play an important role in determining
the spreading or vanishing of the species (Theorem 5.7).
Those theoretical study gives us a method to preserve the rare species who live in the worst environment

Among the work of studying spreading of species, our main interest is that
the domain we investigated is a heterogeneous environment, and the boundary which is governed by a moving boundary
$r=h(t)$. Moreover, we find a threshold number $R^{Fr}_0(D,t)$, which plays the similar importance as the basic reproduction number in epidemiology.
 After the first version of this paper was completed, we have learned a more closely related research in \cite{ZX},
where $b(r)-d(r)$ is defined as $m(x)$ in one space dimension and the authors provided a different way to
understand the dynamics of invasive species by choosing the parameter $D$.
We hope all related work can have useful implications for prediction of biological invasions and preservation of the rare species.


\begin{thebibliography}{99}\setlength{\itemsep}{-1ex}
{\small

\bibitem{AL}
L. J. S. Allen, B. M. Bolker, Y. Lou and A. L. Nevai,
Asymptotic profiles of the steady states for an SIS epidemic reaction-diffusion
model, {\it Discrete Contin. Dyn. Syst. Ser. A.} {\bf 21} (2008), 1-20.

\bibitem{CC}
R. S. Cantrell and C. Cosner, Spatial Ecology via Reaction-Diffusion Equations,
Series in Mathematical and Computational Biology, John Wiley
and Sons Ltd., Chichester, UK, 2003.

\bibitem{CC1}
R. S. Cantrell and C. Cosner, Diffusive logistic equations with indefinite weights:
population models in a disrupted environments, {\it Proc. Roy. Soc. Edinburgh.}
{\bf112A} (1989), 293-318.

\bibitem{CF}
X. F. Chen and A. Friedman, A free boundary problem arising in a
model of wound healing, {\it SIAM J. Math. Anal.} {\bf 32} (2000),
778-800.

\bibitem{Dong}
W. Dong, Positive solutions for logistic type quasilinear elliptic
equations on $\mathbb{R}^n$, {\it J. Math. Anal. Appl.} {\bf290}
(2004), 469-480.

\bibitem{DG}
Y. H. Du and Z. M. Guo, Spreading-vanishing dichotomy in the diffusive
logistic model with a free boundary, II, {\it J. Differential Equations}
{\bf 250} (2011), 4336-4366.

\bibitem{DL}
Y. H. Du and Z. G. Lin, Spreading-vanishing dichotomy in the diffusive
logistic model with a free boundary, {\it SIAM J. Math. Anal.} {\bf
42} (2010), 377-405.

\bibitem{DL2}
Y. H. Du, Z. G. Lin, Erratum: Spreading-vanishing dichotomy in the diffusive
logistic model with a free boundary, {\it SIAM J. Math. Anal.} {\bf
45} (2013), 1995-1996.

 \bibitem{DL1}
Y. H. Du, Z. G. Lin, The diffusive competition model with a free boundary: Invasion of a superior or inferior competitor,
{\it Discrete Contin. Dyn. Syst. Ser. B.}, in press.

\bibitem{DGP}
Y. H. Du, Z. M. Guo and R. Peng, A diffusive logistic model with a free boundary in time-periodic environment,
 {\it J. Funct. Anal.} {\bf 265} (2013), 2089-2142.

\bibitem{DL3}
Y. H. Du and B. D. Lou, Spreading and vanishing in nonlinear diffusion problems
with free boundaries, arXiv:1301.5373, 2013.

\bibitem{DM}
Y. H. Du and L. Ma, Logistic type equations on $\mathbb{R}^N$ by a squeezing  method
involving boundary blow-up solutions, {\it J. London Math. Soc.} (2) {\bf 64} (2001), 107-124.

\bibitem{F}
R. A. Fisher, The wave of advance of advantageous genes, Ann. Eugenics, {\bf 7} (1937), 335-369.

\bibitem{GLL}
H. Gu, Z. G. Lin and B. D. Lou, Different asymptotic spreading speeds induced by advection in
a diffusion problem with free boundaries, Proc. Amer. Math. Soc., in press.

\bibitem{GW}
J. S. Guo, C. H. Wu, On a free boundary problem for a two-species weak competition system,
J. Dynam. Differential Equations (4) {\bf 24} (2012), 873-895.

\bibitem{HH}
W. Huang, M. Han and K. Liu,
Dynamics of an SIS reaction-diffusion epidemic model for disease transmission,
{\it Math Biosci. Eng.} {\bf 7} (2010), 51-66.

\bibitem{KY}
Y. Kaneko and Y. Yamada, A free boundary problem for a reaction-diffusion equation
appearing in ecology, {\it Adv. Math. Sci. Appl.}, {\bf 21} (2011), 467-492.

\bibitem{KPP}
A. N. Kolmogorov, I. G. Petrovsky and N. S. Piskunov, $\grave{E}$tude de l'$\acute{e}$quation de la
diffusion avec croissance de la quantit$\acute{e}$ de mati$\grave{e}$re et son application
$\grave{a}$ un probl$\grave{e}$me biologique, Bull. Univ. Moscou S$\acute{e}$r. Internat.,
Al (1937), 1-26; English transl. in: Dynamics of Curved Fronts, P. Pelc$\acute{e}$(ed.),
Academic Press, 1988, 105-130.

\bibitem{LLW}
C. X. Lei, Z. G. Lin, H. Y. Wang, The free boundary problem describing information
 diffusion in online social networks,
{\it J. Differential Equations} {\bf254} (2013), 1326-1341.

\bibitem{L1}
Y. Lou, Some challenging mathematical problems in evolution of dispersal and population dynamics,
{\it Tutorials Math. Biosci. IV: Evol. Ecol.} (2008), 171-205.

\bibitem{LIN}
Z. G. Lin, A free boundary problem for a predator-prey model, {\it
Nonlinearity}, {\bf 20} (2007), 1883-1892.

\bibitem{LHM}
J. L. Lockwood, M. F. Hoopes and M. P. Marchetti, Invasion Ecology, Blackwell Publishing, 2007.

\bibitem{LSU}
O. A. Ladyzenskaja, V. A. Solonnikov and N. N.
Ural'ceva, Linear and Quasilinear Equations of Parabolic Type, Amer.
Math. Soc, Providence, RI, (1968).

\bibitem{L}
G. Lieberman, Second Order Parabolic Differential Equations, World Sci.,
Singapore (1996) Zbl 0884.35001 MR 1465184.

\bibitem{MCH}
J. Memmott, P. G. Craze, H. M. Harman, P. Syrett and
S. V. Fowler, The effect of propagule size on the invasion of
an alien insect, {\it J. Anim.  Ecol.} {\bf74} (2005), 50-62.

\bibitem{PZ}
R. Peng, X. Q. Zhao, The diffusive logistic model with a free boundary and seasonal succession,
{\it Discrete Contin. Dyn. Syst. A} (5) {\bf 33} (2013), 2007-2031.

\bibitem{PW}
M. H. Protter and H. F. Weinberger, Maximum Principles in Differential Equations,
Springer-Verlag, New York, 1984.

\bibitem{Sa}
D. H. Sattinger, Monotone methods in nonlinear elliptic and parabolic boundary value problems,
{\it Indiana Univ. Math. J.} {\bf21} (1972), 979-1000.

\bibitem{Ske}
J. G. Skellam, Random dispersal in theoretical populations, Biometrika {\bf38} (1951), 196-218.

\bibitem{SK}
N. Shigesada and K. Kawasaki, Biological Invasions: Theory and Practice, Oxford Series in
Ecology and Evolution, Oxford Univ. Press., Oxford, (1997).

\bibitem{Ty}
M. J. Tyler, Australian Frogs A Natural History, Cornell: Cornell Univercity Press, (1994).

\bibitem{Wa}
M. X. Wang, On some free boundary problems of the Lotka-Volterra type prey-predator model, arXiv:1301.2063, 2013.

\bibitem{WZ}
M. X. Wang and J. F. Zhao, A free boundary problem for a predator-prey model with double
free boundaries, arXiv:1312.7751, 2013.

\bibitem{ZX}
P. Zhou and D. M. Xiao, The diffusive logistic model with a free boundary in heterogeneous environment,
{\it J. Differential Equations} {\bf 256} (2014), 1927-1954.

 }
\end{thebibliography}
\end{document}